# Chance Constrained Optimal Power Flow Using the Inner-Outer Approximation Approach

Erfan Mohagheghi, *Member, IEEE*, Abebe Geletu, Nils Bremser, Mansour Alramlawi, *Student Member, IEEE*, Aouss Gabash, *Member, IEEE*, and Pu Li

Department of Simulation and Optimal Processes
Institute of Automation and Systems Engineering
Ilmenau University of Technology
Ilmenau, Germany
erfan.mohagheghi@tu-ilmenau.de, abebe.geletu@tu-ilmenau.de, nils.bremser@tu-ilmenau.de,
mansour.alramlawi@tu-ilmenau.de, aouss.gabash@tu-ilmenau.de, pu.li@tu-ilmenau.de

*Abstract*—In recent years, there has been a huge trend to penetrate renewable energy sources into energy networks. However, these sources introduce uncertain power generation depending on environmental conditions. Therefore, finding 'optimal' and 'feasible' operation strategies is still a big challenge for network operators and thus, an appropriate optimization approach is of utmost importance. In this paper, we formulate the optimal power flow (OPF) with uncertainties as a chance constrained optimization problem. Since uncertainties in the network are usually 'non-Gaussian' distributed random variables, the chance constraints cannot be directly converted to deterministic constraints. Therefore, in this paper we use the recently-developed approach of inner-outer approximation to approximately solve the chance constrained OPF. The effectiveness of the approach is shown using DC OPF incorporating uncertain non-Gaussian distributed wind power.

## Nomenclature

| | |
|---|---|
| $i, j$ | Indices for buses. |
| $B$ | Susceptance. |
| $I(s)$ | Piecewise continuous function. |
| $f$ | Objective function. |
| $\mathbf{g}$ | Equality equations. |
| $\mathbf{h}$ | Inequality equations. |
| $m_1, m_2, \tau$ | Positive constant parameters for analytical approximation. |
| $n$ | Index for samples. |
| $N$ | Total number of samples. |
| $N_b$ | Total number of buses. |
| $OBJ$ | Expected value of objective function. |
| $P$ | Power in a feeder. |
| $P_G$ | Power generation by conventional generators. |
| $P_{G.\min/\max}$ | Lower/upper limit on power generation by conventional generators. |
| $P_L$ | Demand power. |
| $P_{\max}$ | Upper limit on power in a feeder. |
| $Price_G$ | Price of energy generated by conventional generators. |
| $Price_S$ | Price of energy at slack bus. |
| $P_S$ | Power at slack bus. |
| $P_{S.\min/\max}$ | Lower/upper limit on power at slack bus. |
| $P_w$ | Power of a wind farm. |
| $P_{w.F}$ | Forecasted wind power. |
| $S_G$ | Set of conventional generator buses. |
| $T_{CPU}$ | Computation time. |
| $\mathbf{u}$ | Vector of decision variables. |
| $\mathbf{u}_{\min/\max}$ | Lower/upper limits on decision variables. |
| $\mathbf{x}$ | Vector of state variables. |
| $\alpha$ | Probability level. |
| $\beta_w$ | Curtailment factor of a wind farm. |
| $\delta$ | Voltage angle. |
| $\delta_{\min/\max}$ | Lower/upper limit on voltage angle. |
| $\varphi$ | Outer function. |
| $\Theta$ | Parametric function. |
| $\boldsymbol{\xi}$ | Vector of random variables. |
| $\psi$ | Inner function. |
| $\wp$ | Feasible set. |
| $\mathbb{R}^m$ | $m$ dimensional vector space. |

## I. Introduction

Optimal power flow (OPF) [1] has been widely addressed by deterministic approaches which consider the predicted

This work is supported by the Carl-Zeiss-Stiftung.

values of the network variables (e.g., renewable energy generation, demand, prices, etc.) [2, 3]. However, it is not possible to accurately forecast the variables and thus there exist many uncertainties (e.g., demand power [4-6], renewable energy generation [7-13], grid blackouts [14, 15], plug-in electric vehicles [16, 17], etc.) during power system operations. Therefore, network operators have been facing numerous challenges dealing with such uncertainties to ensure not only optimal but also reliable [18, 19] operation strategies.

There are many mathematical models for optimization under uncertainty [20] each of which could be suitable for a specific type of application. For instance, robust optimization and worst-case optimization [21] is frequently used in many applications in which constraint violations are not tolerated. However, in energy networks, there exist some types of constraints (e.g., feeder limits) which are allowed to be violated to some degree and also for a limited time [4].

The promising approach of chance constrained programming [22] is widely used in engineering and finance where uncertainties are common [23-25]. The method was also used to optimize the operation of energy networks under uncertainty[4, 26, 27]. Chance constrained optimization could be used for minimizing the losses and/or maximizing total yield in the network while safeguarding the satisfaction of certain constraints with predefined probability levels. Although formulating the chance constrained optimization problem is advantageous for OPF under uncertainty, it could be in some cases very difficult to solve [28].

For a linear model, if the uncertain variables are normally distributed, there exist deterministic equivalents of chance constraints. Otherwise, there is no direct deterministic representation. Moreover, chance constrained OPF is, in general, a complex problem with uncertain variables described by non-Gaussian probability density function (PDF) [26]. Therefore, the problem should be solved by using an approximation method, e.g., back-mapping [4, 29], sample average approximation (SAA) [30], and inner-outer approximation [31]. Unfortunately, the solution obtained by the SAA method can be infeasible to the chance constrained OPF. On the other hand, back-mapping requires a monotonic property which is commonly not available in power flow problems. Furthermore, to the best of the authors' knowledge, the inner-outer approximation approach has not been utilized in energy networks. The major advantage of this method is that it provides a solution converged from the upper and lower sides, leading to a proof of the feasible solution. Therefore, the main contribution of this paper is using this method to solve the chance constrained OPF under 'non-Gaussian' distributed uncertainties. The results of the stochastic method are then compared to those from the deterministic method confirming the applicability of the approach.

The remainder of the paper is organized as follows. Section II describes chance constrained programming using the inner-outer approximation method. Chance constrained OPF for DC network is formulated in Section III. Section IV presents the results of a case study. The paper is concluded in Section V.

## II. CHANCE CONSTRAINED PROGRAMMING USING THE INNER-OUTER APPROXIMATON METHOD

The OPF problem under uncertainty is generally formulated as follows:

$$\begin{aligned}
(OPF) \quad & \min_{\mathbf{u}} \quad E[f(\mathbf{x},\mathbf{u},\xi)] \\
& s.t. \quad \mathbf{g}(\mathbf{x},\mathbf{u},\xi) = 0, \\
& \quad p(\mathbf{u}) = \Pr\{\mathbf{h}(\mathbf{x},\mathbf{u},\xi) \le 0\} \ge \alpha, \\
& \quad U = \{u \in \mathbb{R}^m \mid \mathbf{u}_{\min} \le \mathbf{u} \le \mathbf{u}_{\max}\},
\end{aligned} \quad (1)$$

where $E[f(\mathbf{x},\mathbf{u},\xi)]$ is the objective function to be minimized with the probabilistic expectation operator $E$, $\mathbf{x}$ is the vector of state variables (e.g., nodal voltages and the power in feeders), $\mathbf{u}$ is the vector of decision variables (e.g., the output power of generators), $\xi$ is the vector of random variables (e.g., renewable energy generations and demands). Since there exist random variables in Eq. (1), the vector of the state variables $\mathbf{x}$ is also random and it could be too expensive to hold the constraints deterministically. Therefore, we use the chance constraint $\Pr\{\mathbf{h}(\mathbf{x},\mathbf{u},\zeta) \le 0\} \ge \alpha$ to satisfy the constraints on state variables by a predefined probability level $\alpha$ with $0.5 \le \alpha \le 1$, where Pr representing a probability measure. Thus, Eq. (1) defines a chance constrained OPF with a feasible set

$$\wp = \{u \in U \mid p(\mathbf{u}) \ge \alpha\}. \quad (2)$$

Since chance constrained optimization problems are generally non-smooth and difficult to solve directly, we approximately solve the chance constrained OPF problem by solving smooth optimization problems. For this, we first define the following function [31]

$$I(s) = \begin{cases} 1, & if \quad s \ge 0 \\ 0, & if \quad s < 0. \end{cases} \quad (3)$$

With this function, we can represent the probability function in Eq. (2) equivalently as

$$\begin{aligned}
\Pr\{\mathbf{h}(\mathbf{x},\mathbf{u},\zeta) \le 0\} \ge \alpha &\equiv E[I(-\mathbf{h}(\mathbf{x},\mathbf{u},\zeta))] \ge \alpha, \\
\Pr\{\mathbf{h}(\mathbf{x},\mathbf{u},\zeta) > 0\} \le 1-\alpha &\equiv E[I(\mathbf{h}(\mathbf{x},\mathbf{u},\zeta))] \le 1-\alpha.
\end{aligned} \quad (4)$$

Note that the function $I(s)$ is not differentiable. Hence, the idea of the inner-outer approximation is to construct a differentiable parametric function $\Theta(\tau, s)$ that resembles the function $I$ and to define an inner approximation (IA) problem

$$\begin{aligned}
(IA) \quad & \min_{u} \quad E[f(\mathbf{x},\mathbf{u},\xi)] \\
& s.t. \quad \mathbf{g}(\mathbf{x},\mathbf{u},\xi) = 0, \\
& \quad \psi(\tau,\mathbf{u}) = E[\Theta(\tau,-\mathbf{h}(\mathbf{x},\mathbf{u},\xi))] \le 1-\alpha
\end{aligned} \quad (5)$$

and an outer approximation (OA) problem

(OA) $\min_{u} E[f(\mathbf{x},\mathbf{u},\boldsymbol{\xi})]$

s.t. $\mathbf{g}(\mathbf{x},\mathbf{u},\boldsymbol{\xi}) = 0,$ (6)

$\varphi(\tau,\mathbf{u}) = E[\Theta(\tau,\mathbf{h}(\mathbf{x},\mathbf{u},\boldsymbol{\xi}))] \geq \alpha$

with the properties
(a) $M(\tau) \subseteq \wp \subseteq S(\tau)$, for $\tau \in (0,1)$,
(b) $\lim_{\tau \searrow 0^+} M(\tau) = \wp = \lim_{\tau \searrow 0^+} S(\tau)$, for $\tau \in (0,1)$,
(c) The problems (IA) and (OA) are differentiable optimization problems,
(d) The solution of (IA) are always feasible to (OPF),
(e) The problem (OA) serves as a parameter tuning strategy for the approximation parameter $\tau$.

From properties (a)-(b), the feasible sets of the approximating problems converge asymptotically to the feasible set of the OPF. Moreover, since the problems (IA) and (OA) are differentiable optimization problems, they can be solved by a gradient-based optimization algorithm. As a result, the cluster points of the solutions of the approximating problems (IA) and (OA) are solutions of the OPF.

In [31] the theoretical foundations of the inner-outer approximation method are given by using the special parametric function

$$\Theta(\tau,s) = \frac{1+m_1\tau}{1+m_2\tau\exp\left(-\dfrac{s}{\tau}\right)},\quad (7)$$

where $\tau$, $m_1$ and $m_2$ are positive constant parameters. In particular, $\lim_{\tau \searrow 0^+} \Theta(s,\tau) = I(s)$. That is, we can approximate the non-smooth function $I(s)$ by the smooth parametric function $\Theta(s,\tau)$, for $\tau \in (0,1)$.

### III. CHANCE CONSTRAINED OPF

In this work, we formulate a DC OPF problem aiming at minimizing the expected cost of power generation by conventional generators (i.e., $P_G(i)$) as well as the expected cost of power imported from a slack bus (i.e., $P_S$). Therefore the objective function is as follows:

$$\min_{\beta_w(i),P_G(i)} E\left[\left(\sum_{i=1}^{N_b} Price_G(i) P_G(i)\right) + Price_S P_S\right],\quad i \in S_G \quad (8)$$

where $Price_G$ is the price of energy generated by conventional generators and $Price_S$ is the price of energy at the slack bus. The objective function in Eq. (8) is subject to the following constraints:

$$\beta_w(i)P_w(i,n) + P_G(i) + P_S(n) - P_L(i,n) =$$
$$\sum_{\substack{i=1\\j=1}}^{N_b} \left(B(i,j)\left(\delta(i,n) - \delta(j,n)\right)\right),\quad n=1,\cdots,N; i \neq j \quad (9)$$

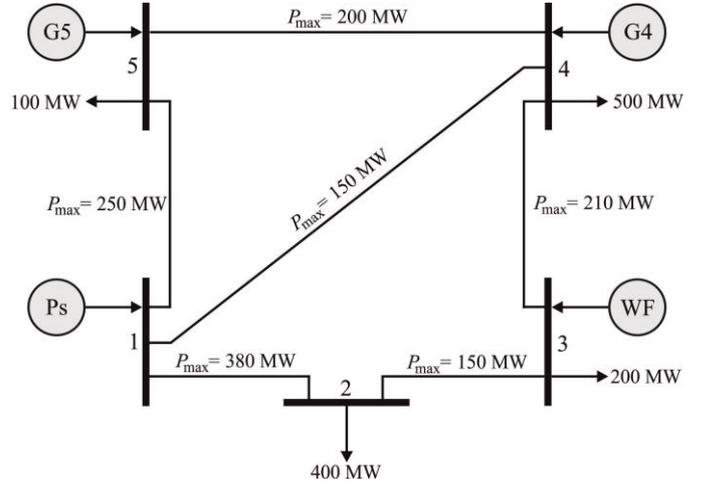

**Fig. 1.** Five-bus network with bus 1 being the slack bus.

where $P_w$ is the power of a wind farm (WF), $\beta_w$ is the curtailment factor of a WF, $P_L$ is demand power, $B$ is the susceptance (i.e., imaginary part of admittance), $\delta$ is the voltage angle, $n$ is the index for samples, $i$ and $j$ are the index of buses. The limitation of power in feeders (i.e., $P(i,j)$) can be expressed as

$$\Pr\{|P(i,j)| \leq P_{\max}(i,j)\} \geq \alpha(i,j),\quad i \neq j \quad (10)$$

where $P_{\max}$ is the upper constraint of the power. The voltage angles are also restricted:

$$\delta_{\min} \leq \delta(n,i) \leq \delta_{\max},\quad n=1,\cdots,N. \quad (11)$$

There is also constraint for power generation of the conventional generators as

$$P_{G.\min}(i) \leq P_G(i) \leq P_{G.\max}(i). \quad (12)$$

The curtailment factor of a WF, i.e., $\beta_w$ is a continuous variable limited to be between 0 and 1:

$$0 \leq \beta_w(i) \leq 1 \quad (13)$$

where $\beta_w = 0$ means full curtailment and $\beta_w = 1$ means no curtailment. Power at the slack bus is limited to

$$P_{S.\min} \leq P_S(n) \leq P_{S.\max}. \quad (14)$$

### IV. CASE STUDY

Fig. 1 shows the network for the case study which is a five-bus network taken and adapted from [32]. Although the inner-outer approximation is a general approach of solving chance constrained optimization problems for linear and nonlinear models, here we solve a DC (linear) OPF under uncertain penetration of wind power. The wind power is described by the Beta PDF [18, 19]. The input data for the optimization is given in Table I and Fig. 2. The objective of the optimization problem is to minimize the total generation costs, while satisfying the feeder power constraints with a probability of at 98%. The problem is coded in GAMS [33] and solved by using CONOPT3 solver.

TABLE I
ENERGY PRICES AND GENERATION LIMITS AT DIFFERENT BUSES

| $Price_S$ ($/MWh) | $Price_G(4)$ ($/MWh) | $Price_G(5)$ ($/MWh) | $Price_w(3)$ ($/MWh) |
|---|---|---|---|
| 15 | 10 | 10 | 0 |
| $P_{S.min}$ (MW) | $P_{G.min}(4)$ (MW) | $P_{G.min}(5)$ (MW) | $P_{w.min}(3)$ (MW) |
| 0 | 0 | 0 | 0 |
| $P_{S.max}$ (MW) | $P_{G.max}(4)$ (MW) | $P_{G.max}(5)$ (MW) | $P_{w.max}(3)$ (MW) |
| 1000 | 400 | 500 | 600 |

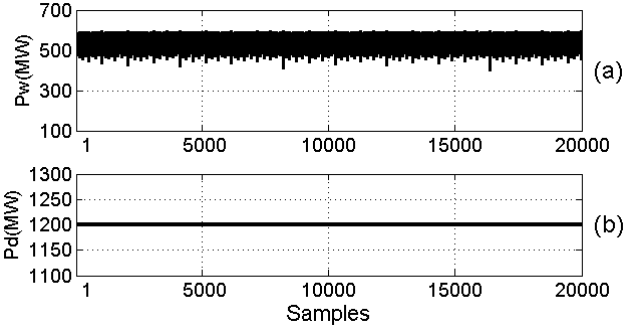

**Fig. 2.** (a) Wind power for 20000 samples. (b) Demand power for 20000 samples.

To show the effectiveness of the method, we carry out optimization for three different cases (i.e., three different forecasted values of wind power): 1) $P_{w.F}(3) = 50$ MW, 2) $P_{w.F}(3) = 300$ MW, and 3) $P_{w.F}(3) = 500$ MW. The resulting curtailment factor of the WF (i.e. $\beta_w$) and the generation of the conventional generators (i.e., $P_G(4)$ and $P_G(5)$) for the three cases are given in Table II. Number of samples, the expected value of objective function and the computation time are given in Table III for the three cases.

The results obtained from the chance constraint OPF is verified using quasi-Monte Carlo sampling method. This gives the true probability of power in the feeders which are given in Tables IV-VI. The trajectories for Case 3 (as a selected case) are shown in Fig. 3 to confirm the effectiveness of the method. The results from our stochastic optimization are compared to those from the deterministic. The significance of our approach can be clearly seen in Fig. 3 where the deterministic approach leads to many violations in feeder constraints.

## V. CONCLUSIONS

Optimal power flow (OPF) is a well-known tool for planning and operation of energy networks. Deterministic approaches have been widely used for OPF in the networks with conventional generation units. However, integration of renewable energies in the networks introduces uncertain generations to the model making those deterministic approaches unsuitable to provide feasible solutions. Therefore, we use the stochastic method of chance constrained optimization to deal with uncertainties associated with wind power. The objective function aims to satisfy predefined levels of constraints satisfaction while minimizing the total costs. However, solving the chance constrained OPF problem is difficult in particular when random parameters are non-Gaussian distributed. To solve this problem, we use the inner-outer approximation method. The effectiveness of the method is confirmed using a linear DC OPF and the advantages are shown over deterministic approaches.

TABLE II
RESULTS OF THE CHANCE CONSTRAINED OPF

|  | $P_{w.F}(3)$ (MW) | $\beta_w(3)$ | $P_G(4)$ (MW) | $P_G(5)$ (MW) |
|---|---|---|---|---|
| **Case 1** | 50 | 1 | 388.2 | 500 |
| **Case 2** | 300 | 0.6986 | 400 | 480.4 |
| **Case 3** | 550 | 0.4495 | 341.78 | 500 |

TABLE III
NUMBER OF SAMPLES, PROBABLITY LEVEL, OBJECTIVE FUNCTION AND COMPUTATION TIME

|  | $N$ | $\alpha$ | OBJ ($) | $T_{CPU}$ (s) |
|---|---|---|---|---|
| **Case 1** | 20000 | 0.98 | 12809.23 | 3.6 |
| **Case 2** | 20000 | 0.98 | 10454.04 | 3.5 |
| **Case 3** | 20000 | 0.98 | 10082.8 | 4.5 |

TABLE IV
THE TRUE PROBABILITY OF POWER IN THE FEEDERS FOR CASE 1

| Bus | 1 | 2 | 3 | 4 | 5 |
|---|---|---|---|---|---|
| 1 | N/A | 0.99385 | N/A | 1 | 0.9999 |
| 2 | 0.99385 | N/A | 1 | N/A | N/A |
| 3 | N/A | 1 | N/A | 0.9789 | N/A |
| 4 | 1 | N/A | 0.9789 | N/A | 1 |
| 5 | 0.9999 | N/A | N/A | 1 | N/A |

TABLE V
THE TRUE PROBABILITY OF POWER IN THE FEEDERS FOR CASE 2

| Bus | 1 | 2 | 3 | 4 | 5 |
|---|---|---|---|---|---|
| 1 | N/A | 1 | N/A | 1 | 0.9783 |
| 2 | 1 | N/A | 0.978 | N/A | N/A |
| 3 | N/A | 0.978 | N/A | 1 | N/A |
| 4 | 1 | N/A | 1 | N/A | 1 |
| 5 | 0.9783 | N/A | N/A | 1 | N/A |

TABLE VI
THE TRUE PROBABILITY OF POWER IN THE FEEDERS FOR CASE 3

| Bus | 1 | 2 | 3 | 4 | 5 |
|---|---|---|---|---|---|
| 1 | N/A | 1 | N/A | 1 | 0.98155 |
| 2 | 1 | N/A | 0.9788 | N/A | N/A |
| 3 | N/A | 0.9788 | N/A | 1 | N/A |
| 4 | 1 | N/A | 1 | N/A | 1 |
| 5 | 0.98155 | N/A | N/A | 1 | N/A |

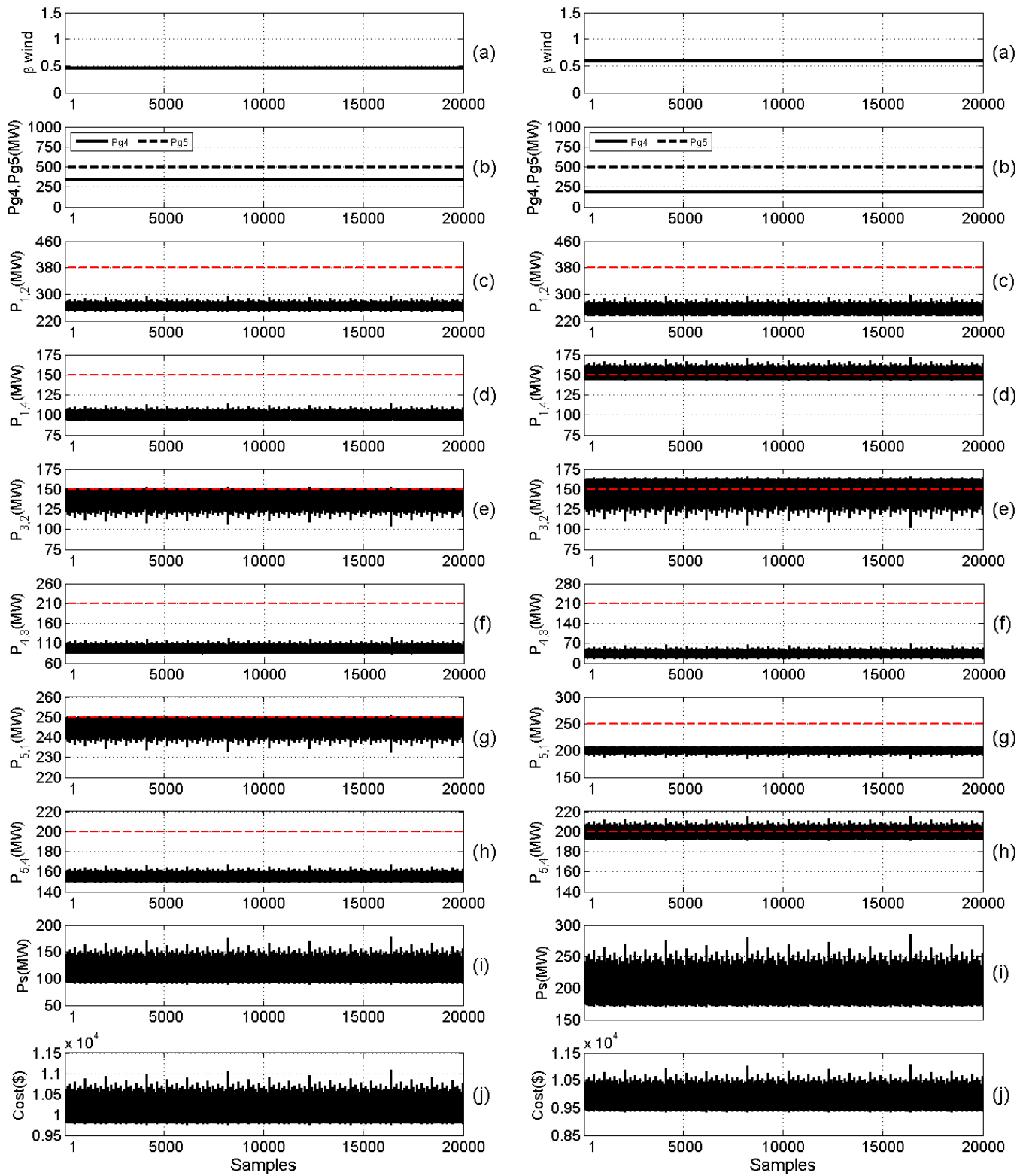

**Fig. 3.** Trajectories for the chance constrained approach (left column) and deterministic approach (right column): (a) Curtailment factors for the wind park. (b) Power generation at buses 4 and 5. (c) Power in the feeder between buses 1 and 2. (d) Power in the feeder between buses 1 and 4. (e) Power in the feeder between buses 3 and 2. (f) Power in the feeder between buses 4 and 3. (g) Power in the feeder between buses 5 and 1. (h) Power in the feeder between buses 5 and 4. (i) Power at the slack bus. (j) Total generation cost.